\documentclass[final]{elsarticle}
\usepackage{latexsym, graphicx, epsfig, amsmath, amsfonts,amssymb,subfigure,algorithm,algorithmicx,booktabs,bm,color,bigints}

\def\Pb{\textrm{Prob}}
\def\mb{\mathbf}

\def\Z{\mathbf{Z}}

\def\I{\mathbb{I}}
\def\i{\mathbf{i}}
\def\j{\mathbf{j}}

\def\t{\tilde}

\def\R{{\mathbb R}}
\def\P{{\mathcal P}}

\def\O{{\Omega}}
\def\o{{\omega}}
\def\A{{\mathcal{A}}}

\def\_#1{{\underline{#1}}}

\newtheorem{remark}{Remark}

\newproof{proof}{Proof}
\begin{document}
\begin{frontmatter}

\title{Efficient failure probability calculation through mesh refinement}

\author{Jing Li}

\address{Pacific Northwest National Laboratory jing.li@pnnl.gov.}

\author{Panos Stinis}

\address{Pacific Northwest National Laboratory panagiotis.stinis@pnnl.gov.}

\begin{abstract}
We present a novel way of accelerating hybrid surrogate methods for the calculation of failure probabilities. The main idea is to use mesh refinement in order to obtain improved local surrogates of low computation cost to simulate on. These improved surrogates can reduce significantly the required number of evaluations of the exact model (which is the usual bottleneck of failure probability calculations). Meanwhile the effort on evaluations of surrogates is dramatically reduced by utilizing low order local surrogates. Numerical results of the application of the proposed approach in several examples of increasing complexity show the robustness, versatility and gain in efficiency of the method.    
\end{abstract}

\end{frontmatter}

\section{Introduction} \label{sec:intro}
Most physical systems are inevitably affected by uncertainties due to natural variabilities or incomplete knowledge about their governing laws. Over the past decade, the engineering research community has realized the importance of advanced stochastic simulation methods for reliability analysis. A quantity of paramount importance in reliability analysis is the calculation of failure probabilities. Mathematically speaking, it is a problem of computing multi-manifold integrals over some failure domains, whose structures are defined by some failure functions, also know as failure modes. The irregular geometry of the failure domain, whose explicit form is usually unavailable, makes the accurate estimation of failure probability very difficult. In applications, the most interesting and relevant cases are usually in high dimensions and with very small failure probability, consequently the problem becomes more challenging. 

A straightforward method for estimating the failure probability (usually denoted by $P_f$) is through Monte Carlo simulations (MC) (cf. \cite{Rubinstein2007}), where a number of samples are drawn and the failure probability is estimated by counting the proportion of the samples lying in the failure domain. The evaluation of the failure function (limit state function) is rather involved and the required sample size used in the evaluations is proportional to $1/P_f$. As a result, the computational cost of brute force MC is usually prohibitively large. To alleviate this inefficiency, various sampling techniques have been explored (e.g., low-discrepancy sampling \cite{Sloan1998,Takhtamyshev2007,Dai2009}, Latin hypercube sampling \cite{Huntington1998,Olsson2003,Pebesma1999}, importance sampling \cite{Au1999,Ibrahim1991,Engelund1993,LiLiXiu_JCP10}, line sampling \cite{Schueller2005,Zio2010}, directional sampling  \cite{Bjerager88,Ditlevsen1990,Nie2000,Nie2004a,Nie2004b}, subset sampling \cite{Au2001,Au2007}). Other alternatives to reduce the simulation effort are non-sampling methods (which are also called approximated methods), such as FORM/SORM(first-order/second-order reliability method) (cf. \cite{Hohenbichler1982-1983,Hohenbichler1987, Fiessler1979, Ditlevsen2007,DerKiureghian1991}). FORM and SORM are based on Talyor expansions (first-order and second-order, respectively) of the failure function around the most probable point on standard Gaussian random space and use asymptotic analytic estimates for the failure probability. Although efficient, approximated methods have reduced accuracy because of the assorted approximating steps adopted in their implementations. 

The response surface method (RSM) is a combination of both simulation methods and approximated methods. To implement RSM, one constructs a surrogate (response surface) which is an approximation of the failure function from deterministic simulations on some design points or heuristic model, and then conducts MC on the obtained surrogate (cf. \cite{Faravelli1989,Rajashekhar1993,Gayton2003,Gupta2004,Paffrath2007,Pulch2010a,jones2012conjunction}). RSM has the advantage of the straightforward implementation of MC and that of the cheaper failure function evaluations on the response surface. However, as was shown in \cite{Li2010}, RMS suffers from robustness issues. To address this problem, a hybrid method was proposed in \cite{Li2010}, where both the response surface and the original underlying failure functions are utilized to construct the approximation of the failure probability. 

The method proposed in this paper combines this hybrid surrogate method with a recently developed mesh refinement method for random space \cite{Li2015,Li2015_2}. 
It is shown in \cite{Li2010} that a better surrogate results in more accurate estimates and requires less computation resources. As a result, to accelerate the hybrid surrogate implementation, constructing a better surrogate should lead to a more efficient algorithm. 

There are two ways that one can use to refine the resolution. One is $p$-refinement which amounts to increasing the degree of the polynomial basis functions used. The second is $h$-refinement which amounts to decreasing the size of the elements used. An accurate high order polynomial surrogate results in a lot of computational effort to be allocated in the evaluation of high order polynomials. However, if we keep the order of the polynomial basis functions only moderately large, we can at the same time refine the random elements into pieces such that the resulting multi-element polynomial expansion has the same error as the high order global surrogate. In the multi-element low order surrogate method, each sample point only needs to be evaluated by a low order polynomial while the multi-element decomposition maintains the accuracy of the surrogate. The adaptive mesh refinement methods proposed in \cite{WanK_JCP05,Li2015,Li2015_2} dynamically decompose the random space according to the evolution of the underlying system. The algorithms automatically allocate more elements near the dynamical important regions which allows to maintain an adequate resolution everywhere. 

With this in mind, we propose two hybrid approaches which fit in the multi-element framework. In the first one we straightforwardly implement the iterative hybrid surrogate approach from \cite{Li2010} and replace the global polynomial surrogate by a multi-element surrogate. The second approach is recommended if one requires very high accuracy of the failure probability estimate. It implements the hybrid surrogate method in each element and then take the sum of the failure probability in each of the elements as the estimate. As the numerical results suggest, the second approach requires more evaluations of the exact failure function but also provides more accurate estimates.    

The rest of the paper is organized as follows. After presenting the formulation of failure probability computation in Section 2, we briefly review
the key ingredients of the present method, including Monte Carlo simulation in Section 3 and the hybrid surrogate approach \cite{Li2010} in Section 4. We rephrase the multi-element gPC expansion and adaptive mesh refinement methods as preliminaries in Section 5. The details of the new method are presented in Section 6, where we explain two versions of our approach. Results of numerical experiments are presented in Section 7 to demonstrate the effectiveness of our algorithms.


\section{Problem formulation} \label{sec:problem}
Let $(\O,\A,\P)$ be a complete probability space, where $\O$ is the event space and $\P$ is the probability measure defined on $\A\in2^{\O}$, the $\sigma-$ algebra of subsets of $\O$ (these subsets are called events). 
Let $\Z = (Z_1,Z_2,\dots,Z_{n_{\Z}}): \O\rightarrow \R^{n_\Z}$ be a $n_{\Z}$ dimensional random vector with joint cumulative distribution function (c.d.f.) $F_{\Z}(z) = \Pb(Z_1\leq z_1,Z_2\leq z_2,\dots,Z_{n_{\Z}}\leq z_{n_{\Z}})$, where $z = (z_1,z_2,\dots,z_{n_\Z})\in \R^{n_{\Z}}$.
A failure domain $\O_{f}$ is usually defined by a limit state function $g$, also called performance function, failure mode or failure constraint as following:
\begin{equation}\label{def:f_domain}
\O_f = \{ \o \in \O : g(\mb{Z(\o)})<0 \}.
\end{equation}
Thus the failure probability can be defined as:
\begin{equation}\label{def:p_f}
P_f  \triangleq \Pb(\o\in \O_f) = \int_{\R^{n_{\Z}}} \I_{\Z(\O_f)}(z)d F_{\Z}(z),
\end{equation}
where $\I$ is the characteristic function satisfying
\begin{equation}\label{def:indicator}
\I_{A} = \left\{ \begin{array} {ccc}1 & \textrm{if} & z\in A,\\ 0 &\textrm{if} & z\notin A.\end{array} \right.  
\end{equation}
Without loss of generality, we focus on the case $\O\subset \R^{n_{\Z}}$ from now on, and let the joint probability density function(p.d.f.) of $\Z$ be $q(z)\triangleq \frac{d F_{\Z}(z)}{d z}$. 

The problem formulation is simple since the calculation of $P_f$ amounts to the calculation of a $n_{\Z}-$dimensional integral. However,  $P_f$ cannot be efficiently calculated by direct numerical integration if $n_{\Z}$ is not small or the geometry of the failure region is not regular. The problem becomes more challenging if the failure probability is small. Since our method is essentially a mixed sampling method, we will briefly introduced the MC method. 

\section{Monte Carlo Simulation}
The Monte Carlo simulation method is widely utilized due to its straightforward implementation and robustness. 
Let $z^{(i)}\in\mathbb{R}^{n_\Z}, i=1,\dots,m$ be a set of samples of the random vector $Z.$ The MC estimate of the failure probability is given by
\begin{equation}\label{failure-prob-monte-carlo-estimate}
P_f^{mc} = \frac{1}{m}\sum_{i=1}^m \I_{\{g(z)<0\}}(z^{(i)}).
\end{equation}
Note that, hereinafter we will, with a slight abuse of notation, use the shorthanded notation
$\{g(z)<0\}$ to stand for the set $\{z:g(z)<0\}$. Although straightforwardly to implement, the MC approach can be costly in practice, since each sample point requires a full-scale simulation of the underlying system. Also, the convergence rate of this estimator is measured by the standard deviation $\sigma^{mc}$ of $P_f^{mc}$, where $\sigma^{mc} = 1/\sqrt{P_f(1-P_f)m}$ \cite{Pradlwarter1999}. Usually a large number of samples is required to obtain an accurate estimate of the failure probability, especially when $P_f$ is small. 

\section{Hybrid surrogate method}\label{sec:hybrid}
A hybrid method was introduced in \cite{Li2010} to combine the robustness of MC with the efficiency of the surrogate approximation. We assume that there exists a surrogate model $\tilde{g}$ that approximates the exact limit state function $g$. The surrogate can be constructed either by numerical approximations or through physicalmathematical reasoning. In general, we assume that the surrogate is an approximation of $g$ with small $L^{p}-$norm error 
\begin{equation}\label{def:lp_bound}
\epsilon_p = \| g(\Z)-\tilde{g}(\Z) \|_{L^p} =( \int| g(z)-\tilde{g}(z) |^{p}q(z) dz)^{1/p}, \quad p\geq 1 .
\end{equation}
In \cite{Li2010} it was shown that straightforward sampling of the surrogate is fundamentally flawed, and may result in erroneous estimates, no matter how accurate the surrogate $\tilde{g}$ is. To alleviate this disadvantage, a hybrid method was proposed, where the direct surrogate sampling result is corrected by results from the limit state function $g$. The key idea of the method is to replace the simulations of the surrogate $\tilde{g}$ when the samples are "close" to the response surface $\tilde{g}=0$ by simulations of the limit state function $g.$ This allows most of the sample simulations of the surrogate $\tilde{g}$ to be kept. 

The hybrid method seeks to estimate the failure probability using the following two integrals,
\begin{equation}\label{def:hybrid_int}
\begin{split}
P_{\gamma}& = \int \I_{\{\tilde{g}<-\gamma\}}(z)q(z)dz,\\
Q_{\gamma}& = \int\I_{\{-\gamma\leq\tilde{g}\leq \gamma\}\cap\{g<0\}}q(z)dz,
\end{split}
\end{equation}
where $\gamma\geq 0$ is a small real number. With these integrals, the estimate of the failure probability by the hybrid method is
\begin{equation}\label{def:hybrid_est}
P_{f}^h = P_{\gamma}+Q_{\gamma}.
\end{equation}
It was proved in \cite{Li2010} that with appropriately chosen threshold parameter $\gamma$, the error of the hybrid estimate $|P_f^h-P_f|$ can be bounded by any prescribed accuracy control threshold $\varepsilon>0.$ The choice of $\gamma$ depends on $\epsilon_p$, the $L^p-$norm ($p\geq 1$) of the error of the surrogate $\tilde{g},$
\begin{equation}\label{def:gamma}
\gamma\geq \frac{1}{\varepsilon^{1/p}}\epsilon_p.
\end{equation}
The hybrid estimate can be computed through MC as follows: Let $\{z^{(i)}\}_{i=1}^M$ be samples drawn from the distribution $q(z)$, then
\begin{equation}\label{pf_h_mc}
\hat{P}_f \triangleq \frac{1}{M}\sum_{i=1}^{M}\left(\I_{\{\tilde{g}<-\gamma\}}(z^{(i)})+\I_{\{|\tilde{g}-g|\leq\gamma\}}(z^{(i)})\cdot\I_{\{g<0\}}(z^{(i)})\right).
\end{equation}
Equation \eqref{pf_h_mc} requires knowledge of the prescribed threshold parameter $\gamma.$ For most practical applications the exact form of $g$ is not available since one only has access to a black box model that provides a result  for each sample. An iterative algorithm letting the algorithm automatically correct the estimate by gradually adding simulations from the exact model served as an alternative of the direct MC computation of the hybrid estimate. In the current work, we only implement the iterative version of the algorithm which avoids the prescription of the threshold parameter $\gamma.$

\vspace{0.5cm}
{\bf Iterative Hybrid Surrogate Algorithm}

Let $m\gg 1$ be the total number of samples, and $S = \{z^{(i)}\}_{i=1}^m$ be the sample set generated according to the distribution $q(z)$. Let $\Delta m,$ called the "step size", be an integer (much) smaller than $m$  and $\eta\geq0$ a small number used at the stopping criterion of the following iteration (also, let integer $\ell\geq0$ be the iteration count):
\begin{itemize}
\item {\bf Initialization:}
\begin{itemize}
\item Set $\ell = 0$, $m_{(\ell)} = 0$,  and $S_{(\ell)} = \emptyset.$
\item Estimate the failure probability using the surrogate model $\tilde{g}.$ That is, for $\ell = 0$, let
\begin{equation}\label{eq:p_h_ini}
\tilde{P}_{(\ell)} = \frac{1}{m}\sum_{i=1}^{m}\I_{\{\tilde{g}<0\}}(z^{(i)}).
\end{equation}
\item Sort $\{|\tilde{g}(z^{(i)})|\}_{i=1}^{m}$ in ascending order.
\end{itemize}
\item {\bf Iteration:} At $\ell-$th iteration ($\ell>0$), do the following:
\begin{itemize}
\item Identify the $m_{(\ell)}+1$ to $(m_{(\ell)}+\Delta m)\wedge m$ elements in the sorted sequence of $|\tilde{g}|$ and their corresponding samples  in the set $S.$ Denote $\Delta S_{(\ell)}$ the set of these samples, and let $S_{(\ell)} = S_{(\ell-1)}\cup \Delta S_{(\ell)}.$
\item Evaluate the exact limit state function $g$ at the sample points in $\Delta S_{(\ell)}.$
\item Update the failure probability estimate using the value of $g$ on $\Delta S_{(\ell)}$,
\begin{equation}\label{eq:p_h_iter}
\tilde{P}_{(\ell)} = \tilde{P}_{(\ell-1)}+\frac{1}{m}\sum_{z\in\Delta S_{(\ell)}}[-\I_{\{\tilde{g}<0\}}(z)+\I_{\{{g}<0\}}(z)].
\end{equation}
\item If $|\tilde{P}_{(\ell)}-\tilde{P}_{(\ell-1)}|\leq \eta$ or $m_{\ell}+\Delta m\geq m,$ 
exit; if not, let $\ell \leftarrow \ell+1$, $m_{(\ell)} \leftarrow m_{(\ell-1)}+\Delta m,$ and repeat the iteration. 
\end{itemize}
\item The failure probability is estimated by
\begin{equation}\label{eq:em_p_h_final}
\tilde{P}_f =\tilde{P}_{(\ell)}.
\end{equation}
\end{itemize}
In \cite{Li2010} it is shown that the accuracy of the surrogate directly impacts the efficiency and accuracy of the estimate of the iterative algorithm. To improve this hybrid algorithm, one approach is to improve the surrogate and another approach is to improve the sampling strategy as presented in \cite{LiLiXiu_JCP10}. If the surrogate is constructed by polynomial expansion, we can use a higher order polynomial as the surrogate. Alternatively, we can divide the random space into elements and use a lower order multi-element polynomial expansion as the surrogate. 

\section{Multi-element gPC surrogate} 
For most problems of practical interest, the limit state function is obtain as the solution to an equation with random parameters, random initial conditions or random boundary conditions. We employ the generalized polynomial chaos (gPC) expansion to construct the surrogate for its exponential convergence rate when the system is smooth with respect to the random parameters. However, it is known that sometimes global polynomial expansion cannot capture the stochastic properties of this solution. In such cases, multi-element polynomial expansion solution will do a much better job. Since in hybrid surrogate method, we need the simulations of most samples on surrogate model, to gain the same accuracy low order multi-element polynomial surrogate is preferred to high order polynomial surrogate. In this section we will introduce the construction of multi-element polynomial surrogate via adaptive gPC mesh refinement. 

\subsection{Generalized polynomial chaos expansion}
The gPC method, an extension of the seminal work on polynomial chaos by R. Ghanem and  (cf. \cite{SpanosG90}), has become one of the most widely used
methods for stochastic computation since its introduction in \cite{XiuK_SISC02}.
Let $\mathbf{i} =(i_1,i_2,\cdots,i_{n_\Z})\in \mathbb{N}_0^{n_\Z}$ be a multi-index with $|\mb{i}|$$ = i_1+i_2+\cdots+i_{n_\Z}$, and let $N\geq 0$ be an integer. The $N$th-degree gPC expansion of $g(\Z)$
takes the following form
\begin{equation}\label{gPC-expansion}
g_N(\Z) = \sum_{|\mb{i}| = 0}^Na_{\mb{i}}\Phi_{\mb{i}}(\Z),
\end{equation}
where $\{a_{\mathbf{i}}\}$ are expansion coefficients that are to be determined, and
$\{\Phi_{\mb{i}}(\Z)\}$ are $n_\Z$-dimensional orthonormal polynomials of degree up to $N$,
satisfying
\begin{equation} \label{gPC-orth}
\mathbb{E}[\Phi_{\bf{i}}(\Z)\Phi_{\bf{j}}(\Z)] = \int_{\Omega}\Phi_{\bf{i}}(z)\Phi_{\bf{j}}(z)q(z)dz
 = \delta_{\bf{i,j}},\qquad 0\leq |\mb{i}|,|\mb{j}|\leq N.
\end{equation}
Here $\delta_{\mb{i,j}}=\prod_{d=1}^{n_\Z}\delta_{i_d,j_d}=1$ if $\mb{i=j}$ and 0 otherwise, and is the
multivariate Kronecker delta function.
The orthogonality relation indicates a correspondence between the distribution of $\Z$ and the type of
the orthogonal polynomial basis $\{\Phi_{\mb{i}}(\Z)\}$. For example, Hermite polynomial chaos corresponds to Gaussian
distribution, Legendre polynomial chaos corresponds to uniform
distribution. For detailed discussions of these relations, see \cite{XiuK_SISC02}.
For a given stochastic system, there are different ways to obtain the expansion coefficients $\{a_{\mathbf{i}}\}_{|\i|=0}^N$ such as stochastic Galerkin method and stochastic collocation method. We will not engage in a detailed discussion on this and refer the interested readers to references such as \cite{Xiu_CICP09}. 

\subsection{Multi-element decomposition in random space}
To deal with irregularities in random space, such as discontinuities, or to obtain a low order polynomial surrogate we discretize the random space $\O\subset \R^{n_{\Z}}$ into a collection of non-overlapping hypercubes. Without loss of generality, let $\Z$ be a $n_{\Z}$-dimensional random vector defined on the probability space $\O$, where $z_i,i=1\cdots,n_{\Z}$ are independent identically distributed (i.i.d) random variables. Let $\O$ be decomposed into $M$ non-overlapping elements as follows:
\begin{eqnarray}\label{decom}
&&B_k = [a_1^k,b_1^k)\times[a_2^k,b_2^k)\times\cdots\times[a_{N_{\Z}}^k,b_{N_{\Z}}^k),\nonumber\\
&&\O = \bigcup_{k=1}^M B_k,\\
&&B_i\cap B_j = \emptyset \quad \text{if } i\neq j,\nonumber
\end{eqnarray}
where $i,j,k = 1,2,\cdots,M$. 

For each random element $B^k$, the local random vector $\Z^k = \Z|_B^k$ has the conditional p.d.f
\[
q^{k}(z^k) = \frac{q(z^k)}{J^k}, \quad k=1,2,\cdots,M, \textrm{ } z^k\in B^k,
\]
where $J^k = \Pb(\I_{B^k}=1) >0$.

Then a multi-element decomposition of the limit state function can be expressed as:
\begin{equation}\label{eq:me_f2}
g(\Z) = \sum_{k=1}^M g(\Z)\I_{B^k}(\Z).
\end{equation}
\subsection{Adaptive gPC mesh refinement} 
Adaptive mesh refinement is a tool to dynamically refine the random space into multi-element according to the evolution or properties of the underlying system.  
From now on in this section, we solve the system in one element ($B^k$) without specially indicating. Assume that the $N$-th order gPC approximation of $g$ can be stated as:
\begin{equation}\label{eq:me_gpc_f}
\tilde{g}(\Z) = \sum_{|\i|=0}^N \tilde{g}_{\i}\Phi_{\i}(\Z).
\end{equation}
Since $\{\Phi_{\i}\}_{|\i|=0}^N$ is a orthonormal basis with respect to the distribution of $\Z$, the local variance is obtained by
\begin{equation}
\sigma^2_{N} = \sum_{|\i|=1}^N(\tilde{g}_{\i})^2.
\end{equation}

If the system is steadied, mesh refinement criterion stated in \cite{WanK_JCP05} can be utilized to guide the refinement in random space. In this method,
a local decay rate of relative error of the approximation in each random element is defined as:
\begin{equation}
\eta = \frac{\sum\limits_{|\i|=N}(\tilde{g}_{\i})^2}{\sigma^2_{N}},
\end{equation}
and the sensitivity of each random dimension is defined as:
\begin{equation}
r_j = \frac{(\tilde{g}{N\mb{e}^j})^2}{\sum\limits_{|\i|=N}(\tilde{g}_{\i})^2},
\end{equation}
where $\mb{e}^j$ is the multi-index such that $\mb{e}^j_i = \delta_{i,j}$, then $N\mb{e}^j$ is the multi-index with $j$-th component $N$. Random element $B^k$ will be split when $$\eta^\alpha \textrm{Prob}(\Z\in B^k) \geq \theta_1,$$ for some $\alpha\in (0,1)$ and error control $\theta_1>0$. If $n_{Z}=1$, $B^k$ will be divided into two equal parts, while if $n_{Z}>1$ all dimensions that satisfy $$r_i\geq \theta_2 \max_{j=1,\dots,n_\Z} r_{j}, \quad i = 1,\dots,n_\Z,$$   
for some $0<\theta_2<1$ will be split into two equal elements.

If the system is dynamical and evolved with respect to time dynamical adaptive mesh refinement algorithm stated in \cite{Li2015,Li2015_2} can be employed to get the multi-element gPC representation. Usually the dynamic of the system is expressed as follows:
\begin{equation}\label{du_t}
u_t(t,z)= \mathcal{L}(t,z;u).
\end{equation}
Assume $\tilde{u} =\sum\limits_{|\i| = 0}^{N} \tilde{u}_{\i}(t)\Phi_{\i}(z)$ is the gPC approximation of the solution to \eqref{du_t} in one element. Then we obtain
\begin{equation}\label{du_t_gPC}
\sum\limits_{|\i| = 0}^{N}\frac{d \tilde{u}_{\i}(t)}{dt}\Phi_{\i}(z) = \mathcal{L}(t,z;\sum\limits_{|\i| = 0}^{N} \tilde{u}_{\i}(t)\Phi_{\i}(z)).
\end{equation}
Taking the Galerkin projection, \eqref{du_t_gPC} becomes
\begin{equation}\label{du_t_galerkin}
\begin{split}
\frac{d \tilde{u}_{\i}(t)}{dt}&= \int \mathcal{L}(t,{z};\sum\limits_{|\j| = 0}^{N} \tilde{u}_{\j}(t)\Phi_{\i}({z}))\Phi_{\i}({z})q({z})dz,\quad |\i| = 0,\dots,N\\
&\triangleq \tilde{R}_\i(\tilde{u}_{\j}(t);|\j|=0,\dots,N),
\end{split}
\end{equation}
here each $\tilde{R}_\i(\tilde{u}_{\j};|\j|=0,\dots,N)$ denotes a function with $\tilde{u}_{\j}$, $|\j|=0,\dots,N$ as its variables, for $|\i|=0,\dots,N$.
If \eqref{du_t_galerkin} has a reliable reduced model of resolvable variables $\hat{u}_{\i}$, $|\i|=0,\dots,N_0$, $N_0<N$ as it is in \cite{Li2015}:
\begin{equation}\label{du_t_reduced}
\frac{d \hat{u}_{\i}(t)}{dt} = \hat{R}_\i(\hat{u}_{\j}(t);|\j|=0,\dots,N_0),\quad |\i|=0,\dots,N_0.
\end{equation}
In this case \eqref{du_t_galerkin} is considered full system while \eqref{du_t_reduced} is the reduced system. We can compute the following quantity
\begin{equation*}
\begin{split}
Q  &= |2\sum_{|\j|=0}^{N_0}\frac{\partial \tilde{u}_{\j}}{\partial t} \tilde{u}_{\j}-2\sum_{|\j|=0}^{N_0}\frac{\partial \hat{u}_{\j}}{\partial t} \hat{u}_{\j}| \\
& = |\sum_{|\i|=0}^{N_0}2\t{R}_\i(\t{u}_{\j}(t);|\j|=0,\dots,N)\t{u}_{\i} - \sum_{|\i|=0}^{N_0}2\hat{R}_\i(\hat{u}_{\j}(t);|\j|=0,\dots,N_0)\hat{u}_{\i} |,
\end{split}
\end{equation*}
which can be considered as "the energy" transferred from the full system to the reduced system, (cf. \cite{Li2015}). The contribution of $i$-th dimension to $Q$ is defined as:
\begin{equation*}
s_i = |2\t{R}_{N_0\mb{e}^i}(\t{u}_{\j}(t);|\j|=0,\dots,N)\t{u}_{N_0\mb{e}^i}-2\hat{R}_{N_0\mb{e}^i}(\hat{u}_{\j}(t);|\j|=0,\dots,N_0)\hat{u}_{N_0\mb{e}^i}|.
\end{equation*}
Random element $B^k$ will be split when $$Q\textrm{Prob}(\Z\in B^k) \geq \theta_1,$$ for some $\theta_1>0$. If $n_{\Z}=1$, $B^k$ will be divided into two equal parts, while if $n_{\Z}>1$ all dimensions that satisfy $$s_i\geq \theta_2 \max_{j=1,\dots,n_\Z} s_{j}, \quad i = 1,\dots,n_\Z,$$   
for some $0<\theta_2<1$ will be split into two equal elements.

If the reduced system \eqref{du_t_reduced} is not available, as long as we have the Galerkin projection \eqref{du_t_galerkin} of the system, we can consider the truncated system up to order $N_0$ as the reduced model (cf. \cite{Li2015_2}) to do the adaptive mesh refinement. $\hat{R}_\i$, $|\i| = 0,\dots,N_0$ is derived by plugging in $\hat{u} =\sum\limits_{|\i| = 0}^{N_0} \hat{u}_{\i}(t)\Phi_{\i}(z)$ in \eqref{du_t} and taking the Galerkin projection. The mesh refinement algorithm is implemented when $Q$ and $s_i$ are calculated with new set of $\{\hat{R}_{\i}\}_{|\i|=0}^{N_0}$.  

After that we have a decomposition of random space and gPC approximation of the limit function (as function of solution to underlying system) in each elements, then the multi-element gPC approximation of limit function is expressed as:
\begin{equation}\label{eq:me_f}
\begin{split}
g(\Z)& = \sum_{k=1}^M g(\Z|_{B^k})\I_{B^k}(\Z)\\
&\approx \sum_{k=1}^M \tilde{g}^k(\Z^k)\I_{B^k}(\Z)=\sum_{k=1}^M \sum_{|\i|=0}^{N}\tilde{g}^k_{\i}\Phi^k_{\i}(\Z^k)\I_{B^k}(\Z).
\end{split}
\end{equation}
\section{Estimate of failure probability in multi-element framework}
Once we get the multi-element expression \eqref{eq:me_f}, the failure probability can be estimated as:
\begin{equation}\label{def:pf_me}
P_f = \int\I_{\O_f}q(z)dz = \sum_{k=1}^M \int\I_{\O_f\cap B^k}q(z)dz.
\end{equation} 
The quantity $$P_f^k \triangleq \int\I_{\O_f\cap B^k}q(z)dz$$ is the contribution to the failure probability from the element $B^k$. 

Assume \eqref{eq:me_gpc_f} is the approximation in element $B^k$, then there exist two hybrid surrogate implementations. The first approach is to consider 
\begin{equation}\label{def:me_f_sur}
\tilde{g} = \sum_{k=1}^M \tilde{g}^k(\Z^k)\I_{B^k}(\Z), 
\end{equation}
as the surrogate and implement the algorithm stated in Section \ref{sec:hybrid}. We call this approach multi-element surrogate global hybrid approach (ME-GHA). Alternatively, since in each random element $B^k$ we have an approximation $\tilde{g}^k,$ the hybrid algorithm can be implemented in each $B^k$ independently and provide an approximation $P_f^k.$ We call this approach multi-element surrogate local hybrid approach (ME-LHA). The ME-LHA algorithm follows 

\vspace{0.5cm}
{\bf Multi-element Surrogate Local Hybrid Iterative Algorithm}

Let $m\gg 1$ be the total number of samples, and $S = \{z^{(i)}\}_{i=1}^m$ be the sample set generated according to the distribution $q(z)$. Let $S^k = S\cap B^k \subset S$ be the samples that are located in the random element $B^k$ and $S^k = \{z^{(k_i)}\}_{i=1}^{m^k}$ for $k=1,\dots,M$, $\sum_{k=1}^M m^k = m.$ Let $\Delta m$ (called the "step size") be an integer (much) smaller than $m$ and $\eta\geq0$ a small number used in the stopping criterion of the following iteration (also let integer $\ell$ be the iteration count):
\begin{itemize}
\item Loop over all elements $B^k$, $k=1,\dots,M$,
\begin{itemize}
\item {\bf Initialization:}
\begin{itemize}
\item Set $\ell = 0$, $m_{(\ell)} = 0$,  and $S_{(\ell)} = \emptyset$.
\item Estimate the failure probability contribution of the $k$-th element using the surrogate model $\tilde{g}^k$. That is, for $\ell = 0$, let
\begin{equation}\label{eq:me_p_h_ini}
\tilde{P}^k_{(\ell)} = \frac{1}{m}\sum_{i=1}^{m^k}\I_{\{\tilde{g}^k<0\}}(z^{(k_i)}).
\end{equation}
\item Sort $\{|\tilde{g}^k(z^{(k_i)})|\}_{i=1}^{m^k}$ in ascending order.
\end{itemize}
\item {\bf Iteration:} At $\ell-$th iteration ($\ell>0$), do the following:
\begin{itemize}
\item Identify the $m_{(\ell)}+1$ to $(m_{\ell}+\Delta m)\wedge m^k$ elements in the sorted sequence of $|\tilde{g}^k|$ and their corresponding samples  in the set $S^k$. Denote $\Delta S_{(\ell)}$ the set of these samples, and let $S_{(\ell)} = S_{(\ell-1)}\cup \Delta S_{(\ell)}$.
\item Evaluate the exact limit state function $g$ at the sample points in $\Delta S_{(\ell)}$.
\item Update the failure probability estimate using the value of $g$ on $\Delta S_{(\ell)}$,
\begin{equation}\label{eq:em_p_h_iter}
\tilde{P}^k_{(\ell)} = \tilde{P}^k_{(\ell-1)}+\frac{1}{m}\sum_{z\in\Delta S_{(\ell)}}[-\I_{\{\tilde{g}^k<0\}}(z)+-\I_{\{{g}<0\}}(z)].
\end{equation}
\item If $|\tilde{P}^k_{(\ell)}-\tilde{P}^k_{(\ell-1)}|\leq \eta$ or $m_{\ell}+\Delta m\geq m^k$, 
\begin{equation}\label{eq:em_p_h}
\tilde{P}_f^k = \tilde{P}^k_{(\ell)},
\end{equation}
exit; if not, let $\ell \leftarrow \ell+1$, $m_{(\ell)} \leftarrow m_{(\ell-1)}+\Delta m$, and repeat the iteration. 
\end{itemize}
\end{itemize}
\item End loop.
\item The failure probability is estimated as:
\begin{equation}\label{eq:em_p_h_final_2}
\tilde{P}_f = \sum_{k=1}^M \tilde{P}_f^k.
\end{equation}
\end{itemize}

\begin{remark}
There are multiple methods to construct the multi-element polynomial surrogates. We can simply divide the random space into elements and construct surrogates in each element as most of multi-element methods do no matter how the system evolves. Also, while we can adaptively refine the random space according to the evolution of the system based on different criteria. For example, one can use the contribution of highest order basis function to the total variance in one element \cite{WanK_JCP05}. Alternatively, one can utilize the rate of change of the activity transferred between full model and reduced model to decide on need for refinement \cite{Li2015} or the rate of change of energy between high order model and low order model in \cite{Li2015_2} and etc. In this paper, we use the latter approach to decide where to place the elements and how to construct the multi-element polynomial surrogates.   
\end{remark}
\begin{remark}
Usually, the surrogates are constructed before implementing the failure probability algorithm. Compared to the sample-sized simulations, we don't take the effort of constructing the surrogate model into account. As long as the surrogate is obtained, the computational cost of the hybrid method is mainly attributed to the simulations on exact model and evaluations of surrogates. Work for evaluation of the polynomial surrogate with the same order is the same. As a result, low order multi-element polynomial surrogates are preferable to high order global polynomials. This means that for the hybrid surrogate method, $h-$refinement is usually better than $p-$refinement.  
\end{remark}
\begin{remark}
For ME-GHA, we directly implement the hybrid method on the multi-element surrogate, while for ME-LHA, the hybrid replacement of surrogate simulations by exact simulations takes place locally in each element. Since at least one iterative step is performed, if we take the same step size as $ME-GHA$, more exact sample simulations are required and thus, better accuracy is achieved. The step size is a prescribed integer, and we can choose a smaller one according to the problem to improve the efficiency. Also, if the total number of exact simulations is limited, an upper bound number of replacements can be set up as stated in \cite{Li2010}.
\end{remark}


\section{Numerical Examples} \label{sec:examples} We present results of the hybrid surrogate approach under the framework of multi-element decomposition for a collection of examples of increasing difficulty: i) a simple example where a global surrogate will not work, ii) a simple linear ODE with a random parameter, iii) the Kraichnan-Orszag three-mode system with random initial conditions and iv) the Burgers equation with random boundary condition.

\subsection{A simple example}
In \cite{Li2010} the authors used counterexamples to demonstrate that direct Monte Carlo simulation of the surrogate model may result in erroneous estimate. One such example is when one considers a step function as the failure function and global gPC expansions as the surrogates. In this subsection, we recalculate the failure probability for this example using a multi-element gPC surrogate. 
Let $Z\sim U[-1,1]$ be a random variable uniformly distributed on $[-1,1].$ Consider the failure function
\begin{equation}\label{ex1_func}
g(Z) = \left\{\begin{array}{rl} -1, & Z\in[-1,0);\\
-0.5, & Z=0;\\
0, & Z\in(0,1]. \end{array}\right.
\end{equation}
By construction $g(Z)\in L^2_{(-1,1)}$. First, we consider the global gPC expansion for $g$ using Legendre polynomials as basis functions. This approximation has the explicit form
\begin{equation}\label{ex1_sur_glo}
\hat{g}^p(Z) = -\frac{1}{2}+\sum_{n=0}^p\frac{(-1)^n(4n+3)(2n)!}{2^{2n+2}(n+1)!(n)!}P_{2n+1}(Z),
\end{equation}
where $P_n$ is the Legendre polynomial of order $n$ (note that $g^p$ is the expansion of order $2p+1$).
 Since $\{P_{n}\}_{n=0}^{\infty}$ is a complete basis for $L^2_{[-1,1]}$, we have 
\begin{equation}\label{ex1_l2}
\lim_{p\rightarrow\infty}\| g-\hat{g}^p \|_{L^2_{[-1,1]}} = 0.
\end{equation}
It is well known that this expansion suffers from Gibbs oscillations. As a result, direct Monte Carlo simulation (MC) on the surrogates cannot give the correct estimate. A numerical simulation with $1,000,000$ samples is implemented and the iterative step size in the hybrid algorithm is chose to be $1,000$ due to the Gibbs phenomenon. The simulation of the same $1,000,000$ samples on the exact model is taken as the reference. 

Results of the direct simulation on surrogates $\hat{g}^p(Z)$ of different $p$ are presented in the second row of Table \ref{tab:ex1}. The results confirm that in this case direct MC on surrogates show erroneous estimates even with very high order(15) global gPC expansion, while with the chosen step size, hybrid method gives the same results as the reference. Number of simulations on the exact model via hybrid method are shown in the third row of Table \ref{tab:ex1}.  Although hybrid method works in this case, half of the samples are have to be simulated by the exact model and thus, it is not efficient. 

Now consider the multi-element surrogate
\begin{equation}\label{ex1_sur_me}
\tilde{g} = \sum_{i=1}^2g_i(Z_i),
\end{equation}
where $Z_1=Z|_{[-1.0)}$, $g_1(Z_1) = -1$, for $Z_1\in[-1,0)$, and $Z_2=Z|_{[0,1]}$, $g_2(Z_2) = 0$, for $Z_2\in[0,1]$.
Since the multi-element polynomial expansion can exactly resolve the original limit state function, applying this surrogate, the hybrid algorithm converges after one iteration, which means only $1,000$ simulations of exact model are required which is $1/1000$ of the sample size compared to $502/1000$ of the sample size when applying hybrid algorithm on global gPC surrogates. The use of a multi-element surrogate for $g$ can improve dramatically the efficiency of the hybrid implementation and at the same time maintain the accuracy of the estimate.

\begin{table}[htbp]
\begin{center}
\begin{tabular}{c|c|c|c}
\hline
$p$ & $0$ & $2$ & $7$\\
\hline
\hline
$\textrm{Prob}(g^p<0)$ & $0.833187$ & $0.773777$ & $0.756490$ \\
$\#$ & $502,000$  & $502,000$  & $502,000$ \\
\hline
\end{tabular}
\caption{Failure probability of step function estimate by global surrogates with different orders}
\label{tab:ex1}
\end{center}
\end{table}

\subsection{A linear one-dimensional ODE}
Consider the linear one-dimensional ODE
\begin{equation}\label{ex:ODE}
\frac{du}{dt} = -Zu,\quad u(0;Z) = u_0,
\end{equation}
where $Z$ is the random input. 
The failure function in this case is defined by $$f(u;Z) = u(T;Z)-u_d.$$ The corresponding failure probability is defined as: $$P_f = \textrm{Prob}[f(u;Z)<0].$$
This equation has the exact solution
\begin{equation}\label{ex:ODE_sol}
u(t,Z) = u_0e^{-Zt}.
\end{equation}
Here we set all parameters the same as in \cite{Li2010} such that $u_0 = 1$, $T = 1$, $u_d = 0.5$ and assume $Z\sim \mathcal{N}({\mu,\sigma^2})$ is a Gaussian random variable with mean $\mu = -2$ and standard deviation $\sigma = 1$. The reference is obtained by direct MC with $1,000,000$ samples.

To construct the multi-element surrogate, we first express the random input as an expansion of a uniform random variable $X$, $X\sim U(-1,1)$. Since $\Phi(x) = \frac{1}{2}+\frac{1}{2}{\rm{erf}}(\frac{x}{\sqrt{2}})$ is the cumulative distribution function (c.d.f.) of standard normal distribution, then 
\begin{equation}\label{uni2nor}
\mu+ \sqrt{2}\sigma{\rm{erf}}^{-1}(X)\sim \mathcal{N}(\mu,\sigma).
\end{equation}
Here ${\rm{erf}}(x) = \frac{\sqrt{2}}{\pi}\int_{0}^x e^{-t^2}dt$ is the error function. Then $Z$ can be approximated by an expansion of Legendre polynomials of a random variable $X\sim U(-1,1)$. We can write 
\begin{equation}\label{nor_apx}
Z\approx \sum\limits_{i=0}^{p}k_{i}L_i(X),
\end{equation}
where $L_i$ is the Legendre polynomial of order $i$. 

We can apply the adaptive mesh refinement algorithm from \cite{Li2015} and we obtain a multi-element solution to \eqref{ex:ODE}. Since our surrogate is a polynomial of a uniform random variable, for comparison our MC samples are taken from uniform distribution and then transformed to a Gaussian random variable by \eqref{uni2nor} to get the reference probability $P_{ref} = 0.003541$. The computational cost for simulation of a sample by multi-element polynomial surrogate of order $p$ is the same as the simulation of a global polynomial surrogate of order $p$. The advantage of our multi-element surrogate lies in employing only low order polynomials to obtain an adequately accurate surrogate. At the same time low order polynomials make the simulations of extremely large sample size affordable.   

Table \ref{tab:ex2} shows the results of the hybrid algorithm with multi-element surrogate and the global surrogate. Given the same sequence of the sample points, each implementation resulted in the same outcome as the reference, however, the number of exact model simulations are different. 
\begin{table}[htbp]
\begin{center}
\begin{tabular}{c|c|c|c|c}
\hline
&$p$ &$3$ & $5$ &$7$\\
\hline
\hline
global surrogate&$\#$  & $105,000$ & $54,700$ & $31,600$ \\
\hline
number of elements &  &$5$ & $5$ & $4$\\
ME-GHA &$\#$ & $3,700$& $3,700$  & $900$ \\
ME-LHA &$\#$  & $4,100$ & $4,100$  & $1200$ \\
\hline
\end{tabular}
\caption{Results of two algorithm implemented on a linear ODE with multi-element polynomials surrogates of different orders p.}
\label{tab:ex2}
\end{center}
\end{table}
The surrogate error has mainly two components. One is the truncation error when we approximate the random input in \eqref{nor_apx}. The second is the numerical error when we evolve the equation \eqref{ex:ODE}. Our adaptive mesh refinement method can reduce the second error. However to reduce the first error we need to increase $p$. To construct the multi-element surrogate, we fixed the tolerance ($TOL_1$) which controls the quality of the mesh refinement. As a result, even with different order $p$, the value of the second error is similar for the two approaches.  

All implementations obtain the same failure probability estimate as the reference. For the computational cost, when $p=3$, the ratio of exact simulations in the multi-element implementation (ME-LHA) to that in the global surrogate implementation is about $1/26$. The ratio is $1/13$ when $p=5$, while $1/26$ when $p=7$. Thus, we see that with the multi-element surrogates, the implementation is much more efficient. As illustrated in the remarks, ME-LHA takes more samples of exact simulations than ME-GHA, and the more elements the larger difference of the number of exact simulations. In this case, ME-LHA and ME-GHA have the same accuracy.

\subsection{The Kraichnan-Orszag three-mode system}
Consider the following system obtained by a linear transformation performed on the original Kraichnan-Orszag three-mode system,
\begin{equation}\label{ex:KO}
\begin{split}
\frac{dy_1}{dt} &= y_1y_3,\\
\frac{dy_2}{dt} &= -y_2y_3,\\
\frac{dy_3}{dt} &= -y_1^2+y_3^2,
\end{split}
\end{equation}
subject to initial conditions
\begin{equation}\label{ex:KO_ini}
y_1(0) = y_1(0;\omega),\quad y_2(0) = y_2(0;\omega),\quad y_3(0) = y_3(0;\omega).
\end{equation}
The solution of the system has a bifurcation depending on the value of the initial conditions $y_1(0)$ and $y_2(0)$. It was shown in \cite{Orszag1967} that a global Wiener-Hermite expansion cannot faithfully represent the dynamics of the system when the random inputs are Gaussian variables. Mesh refinement algorithms \cite{WanK_JCP05,Li2015} can efficiently quantify the uncertainty of the system when the random inputs are uniform random variables whose range includes the bifurcation point. In the current work, we examine only the case of one-dimensional random input. 

We choose the initial conditions
\begin{equation}\label{ex:KO_ini_1d}
y_1(0;\o) = 1,\quad y_2(0;\o) = 0.1\xi(\o),\quad y_3(0;\omega)=0,
\end{equation}
where $\xi\sim U[-1,1]$. The discontinuity (bifurcation) point $y_2=0$ is contained in the random input space. We define the failure function as: 
\begin{equation}\label{ex2:f_func}
f(\o) = y_1(T;\o)-u_d.
\end{equation}
The global gPC expansion cannot accurately describe the solution after about $t=8.$  In our numerical experiments, we fixed $T=15$ and $u_d = 0.03$. The reference result for the failure probability $P_{ref}$ equals $0.102651$ which is obtained by MC with $1,000,000$ samples. The relative error is defined as:
$$
\textrm{relative error} = \frac{|\tilde{P}_f-P_{ref}|}{P_{ref}}.
$$

First we employed the global gPC solutions in the hybrid approach with degree $p=3, 5$, and $7$. We found that the hybrid iterative algorithm cannot obtain a good estimate of the failure probability. In each implementation, the iteration stopped early due to the huge discrepancy between the exact solution and global gPC solution. As a result, the final estimates of the failure probability are not reliable with these global gPC surrogates. 

Next we tested the multi-element surrogates in the hybrid approach. By varying the value of the user-prescribed tolerance $TOL_1$ used to decide the need for refinement, we can obtain different multi-element surrogates. Table \ref{tab:ex3_global} shows results for ME-GHA. Since the tolerance $TOL_1$ controls the error of the surrogate, the stricter the tolerance the better is the surrogate. This means that a smaller number of simulations from the exact model is needed when better accuracy of the estimates is achieved. 

If we fix $TOL_1$, then by increasing $p$, usually we obtained a better surrogate and fewer elements were needed. Table \ref{tab:ex3_local} shows results by ME-LHA.  This means that we implement locally the hybrid method in each random element individually and sum up the contributions to the failure probability from each element to get the estimate of the failure probability. If we take the same step size as in ME-GHA, ME-LHA needs more simulations of the exact model than ME-GHA because at least one iterative step is required in each element.  However these further steps increase the accuracy of the algorithm. Although in each implementation ME-LHA requires more exact simulations than ME-GHA, all implementations achieve the reference estimate.
\begin{table}[htbp]
\begin{center}
\begin{tabular}{c|c|c|c}
\hline
& $TOL_1=10e-3$ & $TOL_1=10e-4$ & $TOL_1=10e-5$\\
\hline
\multicolumn{4}{l}{$p=3$} \\
\hline
No. of elements & $22$ & $38$ & $58$ \\
$\#$ & $6900$ & $500$  &$200$ \\
relative error & $0.06\%$ & $0.16\%$ &0.015\%\\
\hline
\multicolumn{4}{l}{$p=5$} \\
\hline
No. of elements & $12$ & $22$ & $30$ \\
$\#$ & $3900$ & $200$ & $200$ \\
relative error & $0.012\%$& $0.14\%$ & $0.021\%$\\
\hline
\multicolumn{4}{l}{$p=7$} \\
\hline
No. of elements & $10$ & $16$ & $26$ \\
$\#$ & $1400$ & $2200$ & $300$ \\
relative error & $0.33\%$ &$0.038\%$& $0$\\
\hline
\end{tabular}
\caption{Results of ME-GHA for KO system approximated by polynomials with different orders $p$}
\label{tab:ex3_global}
\end{center}
\end{table}

\begin{table}[htbp]
\begin{center}
\begin{tabular}{c|c|c|c}
\hline
& $TOL_1=10e-3$ & $TOL_1=10e-4$ & $TOL_1=10e-5$\\
\hline
\multicolumn{4}{l}{$p=3$} \\
\hline
No. of elements & $22$ & $38$ & $58$ \\
$\#$ & $12245$ & $4700$ & $6029$ \\
\hline
\multicolumn{4}{l}{$p=5$} \\
\hline
No. of elements & $12$ & $22$ & $30$ \\
$\#$ & $3400$ & $3000$ & $3400$ \\
\hline
\multicolumn{4}{l}{$p=7$} \\
\hline
No. of elements & $10$ & $16$ & $26$ \\
$\#$ & $2900$ & $2200$ & $2800$ \\
\hline
\end{tabular}
\caption{Results of ME-LHA for KO system approximated by polynomials with different orders $p$}
\label{tab:ex3_local}
\end{center}
\end{table}


\subsection{Transition Layer for Burgers Equation}
Consider the viscous Burgers equation
\begin{equation}\label{ex:burgers}
\begin{split}
u_t+uu_x &= \nu u_{xx},  \quad x\in[-1,1],\\
u(-1) &= 1+\delta, \quad u(1) = -1.
\end{split}
\end{equation}
Here $\delta$ is a small perturbation of the left boundary ($x = -1$), and $\nu$ is the viscosity. The solution to \eqref{ex:burgers} has a transition layer, defined as the zero of the solution at the steady state. The position of the transition layer is very sensitive with respect to the uncertainty $\delta$ at the left boundary \cite{XiuK_IJNME04,Li2010}. Here we define our failure function as:
\begin{equation}\label{burgers_ffunct}
f(z(\delta)) = -z(\delta)+z_0,
\end{equation}
where $z(\delta)$ satisfies $\lim_{t\rightarrow \infty}u(t,z(\delta);\delta) = 0$, and $\delta\sim U(0,e)$, $e\ll 1$.
The value of $z$ can be found as the solution to a nonlinear system of algebraic equations (\cite{XiuK_IJNME04}).
$$ A\tanh[\frac{A}{2\nu}(1+z)] = 1+\delta,\qquad  A\tanh[\frac{A}{2\nu}(1-z)] = 1. $$
We use the criterion stated in \cite{WanK_JCP05} to adaptively construct the multi-element surrogate. Since we use a collocation method to construct the surrogate in each element, simulations of a new set of collocation points are required whenever refinement takes place (these simulations are counted as simulations of exact model in the cost of implementation). Here we use $21$ collocation points in each element due to the high sensitivity of the transition layer with respect to the random input $\delta$. We fix $z_0 = 0.75$ and calculate the failure probability $P_f = \Pb[f(z)<0]$. Reference failure probability $P_{ref} = 0.127478$ is obtained by MC with $1,000,000$ samples on exact model. We compared the results by global polynomial surrogate of different orders and corresponding multi-element polynomial surrogates in Table \ref{tab:ex4}. While all the algorithms achieve the reference result, the efficiency is much better for the algorithm which uses low order surrogates (p=2,3). 

We compare ME-LHA to the global polynomial surrogate algorithm. We find that the ratio of the number of the exact simulations is $1/40$ when $p=2$, and $1/54$ when $p=3$. With global polynomial surrogate of order $5$, $3921$ simulations of exact model are required which is even larger than the number of exact simulations of multi-element surrogate of order $2$. Evaluation on polynomials of order 2 is much cheaper than that on polynomials of order $5$. 
\begin{table}[htbp]
\begin{center}
\begin{tabular}{c|c|c|c|c|c}
\hline
&$p$ & $2$ & $3$ & $4$ &$5$\\
\hline
\hline
global polynomial&$\#$  & $101,921$ & $62,921$ & $23,421$ & $3,921$ \\
\hline
number of elements & & $9$ & $7$ &$6$ &$5$\\
ME-GHA &$\#$ & $1,757$  & $573$  & $431$ &$389$ \\
ME-LHA &$\#$ & $2,557$  & $1,173$  & $931$ & $799$ \\
\hline
\end{tabular}
\caption{Results of burgers equation via different surrogates (global polynomials and multi-element polynomials) with different orders}
\label{tab:ex4}
\end{center}
\end{table}

\section{Summary}\label{sec:conclusion}
We have presented how the hybrid surrogate method for the calculation of failure probabilities from \cite{Li2010}  can be modified to yield a multi-element random space decomposition. Two versions of the new algorithm were obtained: i) a global hybrid algorithm with multi-element surrogates (ME-GHA) which is a direct implementation of original hybrid surrogate method \cite{Li2010} on multi-element surrogates and ii) a local hybrid algorithm of multi-element surrogates (ME-LHA) in which the hybrid method is implemented on each random element individually. The computational cost for both methods is much smaller than the original hybrid surrogate method. This is achieved firstly by the need for fewer exact model evaluations and secondly by the fact that each surrogate evaluation is cheaper since it uses lower order polynomials. Various numerical examples of increasing difficulty were used to demonstrate the efficacy and robustness of the multi-element hybrid surrogate algorithms.
%
\section*{Acknowledgements} 
This material is based upon work supported by the U.S. Department of Energy Office of Science, Office of Advanced Scientific Computing Research, Applied Mathematics program, Collaboratory on Mathematics for Mesoscopic Modeling of Materials (CM4), under Award Number DE-SC0009280.

\vspace{1.cm}
\bibliographystyle{model1-num-names}
\bibliography{random}

\end{document}